\documentclass[twocolumn]{autart}    

\usepackage{graphicx,bm} 
\usepackage{mathptmx} 
\usepackage{amsmath} 
\usepackage{amssymb}  
\usepackage[mathscr]{euscript}
\usepackage{hyperref}
\usepackage{nicefrac}
\usepackage{times}
\usepackage{enumitem}
\usepackage{picins}
\newtheorem{theorem}{Theorem}  
\newtheorem{lemma}[theorem]{Lemma}  
\newtheorem{proposition}[theorem]{Proposition}    
\newcommand{\subscript}[2]{$#1 _ #2$}

\let\OLDthebibliography\thebibliography
\renewcommand\thebibliography[1]{
  \OLDthebibliography{#1}
  \setlength{\parskip}{0pt}
  \setlength{\itemsep}{0pt plus 0.3ex}
}
\begin{document} 
\begin{frontmatter}

\title{Yield Trajectory Tracking for Hyperbolic\\ Age-Structured Population Systems\thanksref{footnoteinfo}} 
%
\thanks[footnoteinfo]{\textit{Email:} \href{mailto:kevin.schmidt@isys.uni-stuttgart.de}{kevin.schmidt@isys.uni-stuttgart.de} (Kevin Schmidt), \href{mailto:iasonkar@central.ntua.gr}{iasonkar@central.ntua.gr} (Iasson Karafyllis), \href{mailto:krstic@ucsd.edu}{krstic@ucsd.edu} (Miroslav Krstic).}
\author[Stuttgart,SanDiego]{Kevin Schmidt}
\author[Athens]{~~Iasson Karafyllis} 
\author[SanDiego]{~~Miroslav Krstic}
\address[Stuttgart]{Institute for System Dynamics, University of Stuttgart, Waldburgstr. 17/19, 70569 Stuttgart, Germany}                                         
\address[Athens]{Department of Mathematics, National Technical University of Athens, Zografou Campus, 15780, Athens, Greece}             
\address[SanDiego]{Department of Mechanical and Aerospace Engineering, University of California at San Diego, La Jolla, CA 92093-0411, USA} 
\begin{abstract}                          
For population systems modeled by age-structured hyperbolic partial differential equations~(PDEs) that are bilinear in the input and evolve with a positive-valued infinite-dimensional state, global stabilization of constant yield set points was achieved in prior work. Seasonal demands in biotechnological production processes give rise to time-varying yield references. For the proposed control objective aiming at a global attractivity of desired yield trajectories, multiple non-standard features have to be considered: a non-local boundary condition, a PDE state restricted to the positive orthant of the function space and arbitrary restrictive but physically meaningful input constraints.
Moreover, we provide Control Lyapunov Functionals ensuring an exponentially fast attraction of adequate reference trajectories. To achieve this goal, we make use of the relation between first-order hyperbolic PDEs and integral delay equations leading to a decoupling of the input-dependent dynamics and the infinite-dimensional internal one. Furthermore, the dynamic control structure does not necessitate exact knowledge of the model parameters or online measurements of the age-profile. With a Galerkin-based numerical simulation scheme using the key ideas of the Karhunen-Lo{\`e}ve-decomposition, we  demonstrate the controller's performance.\\[-1ex]
\end{abstract}
\end{frontmatter}
\interdisplaylinepenalty=2500
\baselineskip=0.98\normalbaselineskip
\section{Introduction}\label{sec:Intro}
We design an asymptotically tracking control for age-structured chemostats modeled by hyperbolic partial differential equations~(PDEs) with a bilinearly acting input. 
Based on our prior work on the stabilization of constant yield set points, we guarantee global attractivity of output trajectories with an exponential convergence rate. In addition, we developed an efficient numerical scheme based on Galerkin-methods, which ensures accurate asymptotic properties.\\[1ex]
\textbf{Motivation.} In the context of mathematical biology and demography age-structured continuous-time models are a common way of describing the evolution of a certain population with respect to the independent variables of age and time~\cite{Boucekkine.2013,Brauer.2001}. Continuous bioreactors encountered in bioengineering and pharmaceutical research are usually modeled by ordinary differential equations (ODEs)~\cite{Robledo.2012,Smith.1995}. Since multiple ecological concepts like resistance and resilience of ecosystems are closely related to the framework of robustness in system theory, these aspects have been studied rigorously~\cite{Ellermeyer.2001,Karafyllis.2009}. An analysis of the ergodicity problem is given in~\cite{Inaba.1988B,Inaba.1988A}. Moreover, for age-structured models there is an extensive literature on optimal control problems~\cite{Boucekkine.2013,Feichtinger.2003} as well as on the stability of certain PDE models~\cite{Robledo.2012}.\\[1ex]
\textit{The model. }Throughout this paper we focus on the McKendrick-vonFoerster PDE, which is introduced in Section~\ref{sec:Model}. For this setting the dilution rate, which is the ratio of the volumetric flow to the constant volume in the growth chamber, is a natural control variable~\cite{Smith.1995,Toth.2006}. On the other hand, the output is chosen as a weighted integral of the population's age-distribution. As a result, it is possible to represent all products which are proportional to the overall population as well as possibly age-dependent synthesized products. Furthermore, the dependence of the microorganisms' growth rate on the nutrient concentration in the bioreactor is not captured in the model and we hence assume that the nutrient concentration is maintained constant.\\[1ex]
\textit{Time-varying yields. }Having a biomass in mind which is used for the production of antibiotics in pharmaceutical industry, the relevance of the trajectory tracking issue can be elucidated in a demonstrative way. Choosing the yield of an antibiotic as a valid output, its production rate is subjected to external influences like seasonal effects or the current demand. With the prediction of these exogenous factors on an adequate time scale it is possible to determine an optimal production rate governing desired yield trajectories. Due to the fact that a periodic reactor operation may produce a higher yield than the yield achieved by an equilibrium point~\cite{Bittanti.1973}, a special focus is placed on periodical reference trajectories. Beyond this, if the considered bioreactor is a part of a cascaded process, the tracking of predefined trajectories makes it possible to accelerate starting processes and changes of operating points.\\[1ex]
\textbf{Results of the Paper.} 
The present paper extends our prior work~\cite{KaraKrstic.2016}, which focused on the global stabilization of desired equilibrium points of the system class under consideration. We now aim at ensuring the global attractivity of desired yield trajectories and therefore generalize the already established concepts, such that constants set point are included as a special case. For this purpose, a definition of the control objective is given in Section~\ref{sec:ControlTask}. The suggested approach exploits the relation of first-order hyperbolic PDEs to delay models in the sense of integral delay equations (IDEs, see~\cite{Karafyllis.2013}). More specifically, we decompose the PDE problem to an input-dependent finite-dimensional subsystem and an autonomous delay subsystem which is correlated to the microorganisms' reproduction. For special cases of integral kernels we are in the position to construct Control Lyapunov Functionals~(CLFs).\\[1ex]
In contrast to the prior work, we use a two-degrees-of-freedom~(DOF) control structure with a separate feedforward control part evoked by the reference trajectory~\cite{Meurer2009}, as introduced in Section~\ref{sec:ControlDesign}. In this case the feedforward controller does not solely enhance the tracking behavior of the closed loop, but plays an essential role in the overall attractivity concept. The consideration of input constraints is a crucial issue of the present control design assuming a bounded interval for the accessible dilution rate. Moreover, our output-feedback controller does not demand online measurements of the population's entire age-distribution. Even the knowledge of exact system parameters is not necessary, since the controller handles uncertainties in a robust way.
In addition, it is important to guarantee that the PDE state, which represents the population density, remains positive at all times and ages. This fact, in conjunction with a control input directly acting on the whole profile (not simply on the boundary), differentiates our work to other control problems of hyperbolic PDEs~\cite{Bastin.2011,Coron2013}.\\[1ex]
Lastly, we provide simulation results of the closed-loop system in Section~\ref{sec:Simulation} with a Galerkin-based simulation scheme, which conserves important system properties even at low orders and enables independent age and time discretizations.\\[2ex]
%
\textbf{Notation. }\\[-4ex]
\begin{itemize}
  \item The set $\mathbb{R}^+$ denotes all positive-valued real numbers, $\mathbb{R}^+_0$ all non-negative real numbers.\\[-1.5ex]
	\item The inner product of $L^2$ is denoted~$\langle f, g\rangle:=\int_0^A f(a)g(a)\text{d}{a}$ where $f,g\in L^2([0,A])$.\\[-1.5ex]
	\item $\left\|f(a)\right\|_{\infty}=\max_{a\in[0,A]} | f(a) |$ is the maximum- resp. $L^{\infty}$-norm for $f\in \mathcal{C}^0([0,A])$.\\[-1.5ex]
	\item $\mathcal{K}_\infty$ is the class of all strictly increasing, unbounded functions~$\kappa\in\mathcal{C}^0(\mathbb{R}^+_0;\mathbb{R}^+_0)$ with $\kappa(0)=0$.
\end{itemize}
\begin{itemize}
	\item The saturation function with respect to~$f\in[D_\text{min},D_\text{max}]$ is defined $\text{sat}(f)=\min(D_\text{max},\max(D_\text{min},f))$; other intervals are explicitly denoted as an index.\\[-1.5ex]
	\item Given the functions $f:~\mathbb{R}^+\times X\to\mathbb{R}$ , $z:~\mathbb{R}^+\to X$ with the metric space~$X$, we define the right temporal Dini-derivative
$	\dot{f}^+(t,z(t)):=\overline{\lim}_{h\to0^+} \frac{f(t+h,z(t+h))-f(t,z(t))}{h}$.\\[-1.5ex]
\item For any~$S\subseteq\mathbb{R}$ and $A>0$, $P\mathcal{C}^1([0,A];S)$  denotes the class of all functions~$f(a)\in\mathcal{C}^0([0,A];S)$ for which there exists a finite (or empty) set $B\subseteq(0, A)$ such that: (i) the derivative $f'(a)$ exists at every $a\in(0, A)\setminus B$ and is a continuous function on $(0, A) \setminus B$ , (ii) all meaningful right and left limits of $f'(a)$ when $a$ tends to a point in $B\cup\{0, A\}$ exist and are finite.
\end{itemize}
%
\section{Age-Structured Population Models} \label{sec:Model}
Consider the McKendrick-vonFoester PDE~\eqref{PDE} valid in the age-time domain~$(a,t)\in(0,A)\times(0,\infty)$
\begin{align}
 \label{PDE}\frac{\partial x(a,t)}{\partial t}+\frac{\partial x(a,t)}{\partial a} 	&=-[\mu(a)+D(t)]x(a,t)\\
 \label{BC}x(0,t)																&=\int_0^A k(a)x(a,t)\text{d}{a}\\
 \label{IC}x(a,0) 																&= x_0(a) \end{align}
which describes the evolution of the population density $x: ~[0,A]\times[0,\infty)\to\mathbb{R}^+$ as a part of an initial-boundary value problem (IBVP) on the same domain with an arbitrary large but finite maximum reproductive age $A>0$. Strictly speaking, the state~$(x[t])(a)=x(a,t)$ describes the density of the overall population which has reached a specific age~$a$ at a certain time~$t$. In addition, the function $\mu(a)$ denotes the age-dependent mortality rate and $D(t)$ the dilution rate which is the control input.
In particular, the non-local boundary condition~(BC)~\eqref{BC} is valid for $t\geq0$ and models the production of new-born individual $v(t)=x(0,t)$ determined by the birth modulus resp. the kernel~$k(a)$. Furthermore, \eqref{IC} is the initial condition, i.e. the initial distribution of the population density in the age-domain~$[0,A]$ at $t=0$. In addition, the output is defined by the equation
\begin{equation}
\addtolength{\abovedisplayskip}{-2mm}
\addtolength{\belowdisplayskip}{-2mm}
 \label{OE}y(t)=\int_0^A p(a)x(a,t)\text{d}{a},
\end{equation}
which takes the possibly age-dependent production rate~$y(t)$ of a specific (bio)chemical species into account. For instance, the choice~$p(a)=1$ includes the overall population as a valid output.\\
The distributed parameter system~$\Sigma_x$ given by~\eqref{PDE}--\eqref{OE} with the input~$D(t)$ and the output $y(t)$ is of bilinear single-input-single-output type. Subsequently, we introduce three assumptions in order to guarantee the existence of a meaningful unique solution of the IBVP~\eqref{PDE}--\eqref{BC}  aware of the state and input constraints~(see also~\cite{KaraKrstic.2016}):
\begin{enumerate}[label=(\subscript{$A$}{\arabic*})]
    \item The parameters functions are restricted to $k, p\in\mathcal{P}$ and~$\mu\in\mathcal{C}^0([0,A];\mathbb{R}^+_0)$, where $\mathcal{P}:=\{f\in\mathcal{C}^0([0,A];\mathbb{R}^+_0)\big|$$~\langle 1, f \rangle> 0\}$.\\[-1ex]
\end{enumerate}
\begin{enumerate}[label=(\subscript{$A$}{\arabic*})]
\setcounter{enumi}{1}
   \item The control~$D(t)$ takes values in $[D_{\min},D_{\max}]$ $\subset\mathbb{R}^+_0,$ where $D_{\min}<D_{\max}$.\\[-1ex]
	\item The initial condition (IC)~\eqref{IC} is compatible with~\eqref{BC}, i.e.
$x_0\in\mathcal{X}:=\{f\in P\mathcal{C}^1([0,A];\mathbb{R}^+)\big|$ $~f(0)=\langle k, f \rangle> 0\}$.
\end{enumerate}
%
\section{Control Objective and Delay Model} \label{sec:ControlTask}
The asymptotic tracking of a reference trajectory~$y_\text{ref}(t)$ with respect to the output~$y(t)$ given by~\eqref{OE} defines the key objective of the contribution. For designing an asymptotic tracking control in the context of system~$\Sigma_x$, consider the 2-DOF control structure satisfying assumption~(A$_2$) a priori\\[-3ex]
\begin{figure}
\centering
\includegraphics[width=0.48\textwidth]{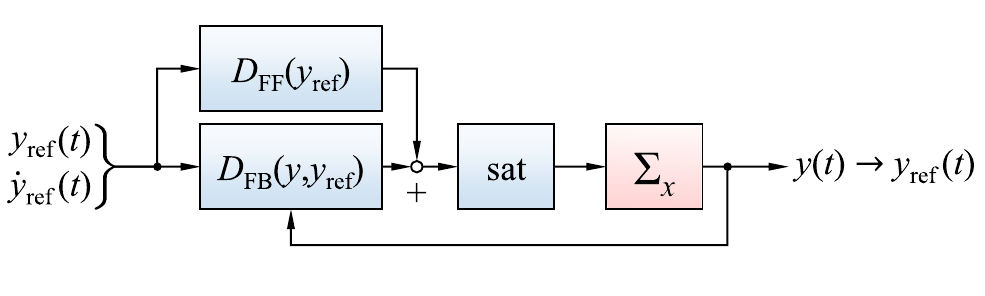}
\caption{2-DOF control structure for the plant~$\Sigma_x$ \eqref{PDE}--\eqref{OE} with a feedforward component $D_\text{FF}(y_\text{ref})$ and a feedback part $D_\text{FB}(y,y_\text{ref})$.} 
\label{fig:2DOF_ControlStructure}                                 
\end{figure}
\begin{align}
\label{2DOF_controller}
D\left(y,y_\text{ref}\right) =  \text{sat}\Big(D_\text{FF}(y_\text{ref})+D_\text{FB}(y,y_\text{ref})\Big)\\[-4ex]\nonumber
\end{align}
consisting of a feedforward~$D_\text{FF}(y_\text{ref})$ and a feedback component~$D_\text{FB}(y,y_\text{ref})$ as sketched in Fig.~\ref{fig:2DOF_ControlStructure}. More specifically, the desired global asymptotic attractivity w.r.t.~$y_\text{ref}(t)$ for all valid initial conditions is defined as follows\\[-3ex]
\begin{align}
\lim_{t\to\infty}|y(t)-y_\text{ref}(t)| = 0~~~~~\forall x_0\in\mathcal{X}.\\[-4ex]\nonumber
\end{align}
There are basically two different types of relevant reference trajectories in the context of chemostats. The first type focuses on the planning of a trajectory which reaches a set point $y_\Delta$ in a desired transition time~$t_\Delta$ starting from an arbitrary initial profile and is motivated by starting processes or set point changes during an operating chemostat. In addition the stabilization of periodical trajectories is relevant due to an efficient harvesting of synthesized products~\cite{Mazenc.2008} and is considered as a second case.\\[1ex]
Before starting with the control design, we analyze the equilibrium profiles and introduce a delay system to describe the dynamics of system \eqref{PDE}--\eqref{OE}. Therefore, define the equilibrium dilution rate~$D^*$ as the unique solution of\\[-3ex]
\begin{align}
\label{LotkaSharpe}
1=\int_0^A k(a) \e^{-D^*a-\int_0^a\mu(\alpha)\text{d}{\alpha}}\text{d}{a} =: \int_0^A  \tilde{k}(a) \text{d}a .
\end{align}
It follows from a steady state analysis of~$\Sigma_x$:\,\eqref{PDE}--\eqref{OE} that the necessary condition $D(t)=D^*$ corresponds to a continuum of equilibrium profiles~\cite{KaraKrstic.2016}, which are proportional to
\begin{align}
\label{x*(a)}
x^*(a) = \frac{\e^{-D^*a-\int_0^a \mu(\alpha)\text{d}{\alpha}}}{\int_0^A p(s)  \e^{-D^*s-\int_0^s \mu(\alpha)\text{d}{\alpha}}\text{d}{s}} .
\end{align}
As a prerequisite for the following proposition, we use this fact to define the kernel $g(a)$ and a reference profile of the PDE state~$x_\text{ref}(a,t)$ with the property $y_\text{ref}(t)=\langle p,x_\text{ref}[t]\rangle$:
\begin{align}
\label{g(a)}
 g(a) &=x^*(a)\cdot p(a),&&\text{for all }a\in[0,A]\\
x_\text{ref}(a,t) &=x^*(a)\cdot{y_\text{ref}(t)},&&\text{for all }a\in[0,A], t\geq0.
\end{align} 
In the following proposition we refer to the open-loop system, i.e. it is valid for any piecewise continuous input $D(t)$ satisfying assumption~(A$_2$). Its proof is omitted, since it follows the arguments of the proof of Lemma~\ref{lem:delay_cl}, which addresses the controlled system (see Section~\ref{sec:ProofMain})
\begin{proposition}[delay model]\label{thm:Trafo} Define the functional
\begin{align}
\label{def_Pi}
\Pi(f):=\frac{\left\langle \pi, f \right\rangle}{\left\langle \pi, x^* \right\rangle}=\frac{\int_0^A \pi(a) f(a)\emph{d}a}{\int_0^A \pi(a) x^*(a)\emph{d}a}
\end{align}
%
where $\pi:~[0,a]\to\mathbb{R}^+$ is the continuous function
\begin{align}
\label{pi_0}
\pi(a) :=  \int_a^A k(s) \e^{D^*(a-s)+\int_s^a \mu(\alpha)\emph{d}{\alpha}}\emph{d}{s}.
\end{align}
Given an IC $x_0\in\mathcal{X}$ for the PDE system~\eqref{PDE}--\eqref{OE}, we define
\begin{align}
\label{IC_etapsi}
\eta_0 = \ln\,\frac{\Pi(x_0)}{y_\emph{ref}(0)},\,
\psi_0(a) = \frac{x_0(a)}{x^*(a)\Pi(x_0)}-1,\,a\in[0,A].
\end{align}
For every piecewise continuous input $D:~\mathbb{R}^+_0\to[D_\emph{min},D_\emph{max}]$, consider the solution of the delay model
\begin{align}
\label{ODE_eta}\dot{\eta}(t) &=D^*-\frac{\dot{y}_\emph{ref}(t)}{y_\emph{ref}(t)}-D(t)\\
\label{IDE_psi}\psi(t) 			 &= \int_0^A \tilde{k}(a) \psi(t-a)\emph{d}a.
\end{align}
with the ICs $\eta(0)=\eta_0$ and $\psi(-a)=\psi_0(a)$ for $a\in[0,A]$. Define the functions for $(a,t)\in[0,A]\times\mathbb{R}^+_0$,
\begin{align}
\label{x_psieta}
\Phi_x(a,t) &= x_\emph{ref}(a,t) \e^{\eta(t)} \left[1+\psi(t-a)\right]\\
\label{y_psieta}\Phi_y(t) &=y_\emph{ref}(t) \e^{\eta(t)} \left[1+\int_0^Ag(a)\psi(t-a)\emph{d}a\right].
\end{align}
Then, the unique solution of the PDE model~\eqref{PDE}--\eqref{OE} corresponding to the input~$D$ and the IC~$x_0\in\mathcal{X}$ is given by $x(a,t)=\Phi_x(a,t)$ for $(a,t)\in[0,A]\times[0,\infty)$. Moreover, the output~$y(t)$ is given by~$y(t)=\Phi_y(t)$ for $t\geq0$.\\[-3ex]
\end{proposition}
With~\eqref{IC_etapsi}--\eqref{IDE_psi} we decompose the PDE dynamics to an input-dependent ODE-state~$\eta(t)$ and an infinite-dimensional internal coordinate~$\psi(t-a)$ which cannot be affected by the input~$D(t)$. It can be seen from~\eqref{y_psieta}, that a steady state of the delay model $\eta(t)=\psi(t-a)\equiv0$, implies an exact output tracking, i.e. $y(t)=y_\text{ref}(t)$. With the definition\\[-3ex]
\begin{align}
\label{def_delta}
\delta(t):=\ln\left(1+\int_0^A g(a)\psi(t-a)\text{d}a\right)
\end{align}
we rewrite~\eqref{y_psieta} and get following result for the logarithmic tracking error as the basis for the control design
\begin{align}
\label{Y_psidelta}
\ln\left(\frac{y(t)}{y_\text{ref}(t)} \right)=\eta(t)+\delta(t).
\end{align}
%
%
%
\section{Nonlinear Controller Design} \label{sec:ControlDesign}
For the design of an appropriate controller, we state a second set of assumptions:
\begin{enumerate}[label=(\subscript{$B$}{\arabic*})]
    \item The equilibrium dilution rate~$D^*$ defined by~\eqref{LotkaSharpe} satisfies~$D^*\in(D_\text{min},D_\text{max})$.\\[-1ex]
		\item Assume that $y_\text{ref}$ belongs to the class of valid reference trajectories $\mathcal{Y}$, which is defined as all positive-valued and continuously differentiable functions $\mathcal{C}^1(\mathbb{R}^+_0;\mathbb{R}^+)$ with the additional property  
		\begin{equation}
\label{ValidYref}
 D^*-D_{\max}<\inf_{t\geq 0} \frac{\dot{y}_\text{ref}(t)}{y_\text{ref}(t)}\leq\sup_{t\geq 0} \frac{\dot{y}_\text{ref}(t)}{y_\text{ref}(t)} < D^*-D_{\min}.
\end{equation}
	  \item We assume the existence of a constant $\lambda>0$ which fulfills the inequality with $\tilde{k}(a)$ defined by~\eqref{LotkaSharpe}
 \begin{equation}
\int_0^A\left|\tilde{k}(a)-\lambda \frac{\int_a^A \tilde{k}(s)\text{d}{s}}{\int_0^A s \tilde{k}(s)\text{d}{s}}\right|\text{d}{a}<1.
\end{equation}
	\item The system parameters $\mu(a),k(a),p(a)$ as well as the initial profile~$x_0(a)$ are assumed to be unknown. As a result of~\eqref{LotkaSharpe}, the equilibrium input~$D^*$ is also unknown.\\[-3ex]
	\end{enumerate}
While assumption (B$_1$) ensures the reachability of the plant's steady states, the extra condition~\eqref{ValidYref} involved in assumption (B$_2$) is needed to guarantee the asymptotic tracking in the presence of input constraints. Note that if the birth modulus~$k\in\mathcal{P}$ has a finite number of zeros in the age-domain~$a\in[0,A]$, assumption~(B$_3$) holds true by virtue of Proposition 2.3 in~\cite{KaraKrstic.2016}, no matter what $\mu(a)$ and $D^*$ are. Moreover, this assumption is not necessary to prove stability but to explicitly construct CLFs. \\[1ex]
Consider the nonlinear dynamic output controller specifying the 2-DOF structure in~\eqref{2DOF_controller}  with the control gain~$\gamma>0$ and the observer gains $l_1,l_2>0$ for $t>0$ 
\begin{align}
\label{D_FF}
D_\text{FF}(t) &= -\frac{\dot{y}_\text{ref}(t)}{y_\text{ref}(t)}\\
\label{D_FB}
D_\text{FB}(t) &=z_2(t)+\gamma \ln\frac{{y}(t)}{y_\text{ref}(t)}\\
\label{z_dot}
\dot{\bm{z}}(t)&=\underbrace{\begin{bmatrix} -l_1 & 1 \\ -l_2 &0 \end{bmatrix}}_{\bm{L}}\bm{z}(t)+\begin{bmatrix} l_1\ln\frac{{y}(t)}{y_\text{ref}(t)}+D_\text{FF}(t)-D(t) \\  l_2\ln\frac{{y}(t)}{y_\text{ref}(t)}\end{bmatrix}\\
\label{D_CL}
D(t) &= \text{sat}\left(-\frac{\dot{y}_\text{ref}(t)}{y_\text{ref}(t)} +z_2(t)+\gamma \ln\frac{{y}(t)}{y_\text{ref}(t)}\right).
\end{align}
The proposed controller~\eqref{D_FF}--\eqref{D_CL} consists of three components: Firstly, the feedforward control $D_\text{FF}(t)$ is constructed by means of the scalar dynamics~\eqref{ODE_eta}. On the other hand, the feedback part~\eqref{D_FB} consists of a proportional control of the logarithmic error~\eqref{Y_psidelta} and an observer-ODE~\eqref{z_dot} which adapts the controller to the unknown equilibrium dilution rate~$D^*$.
Hereby, the observer state $\bm{z}(t)\in\mathbb{R}^2$ has the IC~$\bm{z}_0\in\mathbb{R}^2$ and its second component $z_{02}\in\mathbb{R}$ is the initial guess for the equilibrium dilution rate $D^*$. For $y(t)=y_\text{ref}(t)$ and $D(t)=D^*+D_\text{FF}(t)$, $\bm{z}(t)$ has a steady state $\bm{z}^*:=[0,D^*]^{\intercal}$. As a result we define the observer error
\begin{align}
\label{def_e(t)}
\bm{e}(t)=\bm{z}(t)-\Big[\ln\frac{\Pi(x[t])}{y_\text{ref}(t)},~D^*\Big]^{\intercal}\hspace{-1mm}=\bm{z}(t)-[\eta(t),~D^*]^{\intercal}
\end{align}
%
%
The control law~\eqref{D_FF}--\eqref{D_CL} guarantees an attractivity with exponential convergence of an admissible reference trajectory $y_\text{ref}\in\mathcal{Y}$ as characterized by the subsequent theorem:\\[-3ex]
\begin{theorem}[robust global attractivity]\label{thm:GES}
Consider the system $\Sigma_x:$~\eqref{PDE}--\eqref{OE} under assumptions~$(A_1)$--$(A_3)$ and~$(B_1)$--$(B_4)$. For every~$y_\emph{ref}\in\mathcal{Y}$, there exists a constant~$L(y_\emph{ref})>0$ and a function~$\kappa\in\mathcal{K}_\infty$~(independent of $y_\emph{ref}$), such that the unique solution of the closed loop~\eqref{PDE}--\eqref{OE} with the controller~\eqref{D_FF}--\eqref{D_CL} satisfies the following estimate for any $x_0\in\mathcal{X}$ and $\bm{z}_0\in\mathbb{R}^2$ and all~$t\geq0$:\\[-3ex]
\begin{align}
\label{thm:Linf}
\left\|\ln\frac{x(a,t)}{x_\emph{ref}(a,t)}\right\|_{\infty}+\left|\bm{e}(t)\right|\leq& \kappa\left(\left\|\ln\frac{x_0(a)}{x_\emph{ref}(a,0)}\right\|_{\infty}+\left|\bm{e}_0\right|\right) \nonumber\\ &\times \emph{e}^{-\frac{L(y_\emph{ref})}{4}t} 
\end{align}
Moreover, for any constants $l_1,l_2,p_1,p_2>0$ satisfying
\begin{align}
\label{ineq_p12}
(2+l_1p_1-2l_2p_2)^2 < 8l_1p_1-4l_2p_1^2,~~~~
p_1^2 < 4p_2,\\[-4ex]\nonumber
\end{align}
there exists a sufficiently small constant $\sigma>0$, sufficiently large constants $\alpha_1,M>0$, such that for every $\alpha_2>0$ the functional
$~V:~\mathbb{R}^2\times\mathcal{X}\times\mathbb{R}^+_0\to\mathbb{R}^+_0$
defined by\\[-1ex]
\begin{equation}
\label{CLF_V}
V(\bm{z},x,t)=\left(\ln\frac{\Pi(x)}{y_\emph{ref}(t)}\right)^2+\alpha_1 \sqrt{Q(\bm{z},x,t)}+\alpha_2 Q(\bm{z},x,t)
\end{equation}
with\\[-4ex]
\begin{align}
Q(\bm{z},x,t)&=	\begin{bmatrix}z_1-\ln\frac{\Pi(x)}{y_\emph{ref}(t)}\\ z_2-D^*\end{bmatrix}^{\intercal}\begin{bmatrix}1 & -\frac{p_1}{2} \\ -\frac{p_1}{2} & p_2\end{bmatrix}\begin{bmatrix}z_1-\ln\frac{\Pi(x)}{y_\emph{ref}(t)}\\ z_2-D^*\end{bmatrix}
 \nonumber\\[1ex]
						&+\frac{M}{2}\left(\frac{\left\|\e^{-\sigma a}\left(\frac{x(a)}{x_\emph{ref}(a,t)}-\frac{\Pi(x)}{y_\emph{ref}(t)}\right)\right\|_{\infty}}{\min \left\{ \frac{\Pi(x)}{y_\emph{ref}(t)},\min_{a\in[0,A]} \frac{x(a)}{x_\emph{ref}(a,t)}\right\}}\right)^2 
						\end{align}
is a Lyapunov functional for the closed loop system~\eqref{PDE}--\eqref{OE},~\eqref{D_FF}--\eqref{D_CL}. More specifically, the differential inequality
\begin{equation}
\label{thm:dotV}
\dot{V}^+(\bm{z}(t),x[t],t)\leq -L(y_\emph{ref})\cdot\frac{V(\bm{z}(t),x[t],t)}{1+\sqrt{V(\bm{z}(t),x[t],t)}}
\end{equation}
holds for all $t\geq 0$ and every solution $(x[t],\bm{z}(t))\in\mathcal{X}\times\mathbb{R}^2$ of the closed loop~\eqref{PDE}--\eqref{OE} with the controller \eqref{D_FF}--\eqref{D_CL}.
\end{theorem}
Theorem~\ref{thm:GES} does not only guarantee the global asymptotic attractivity of the output $y(t)\to y_\text{ref}(t)$, but it determines the whole PDE state $x(a,t)$ of the closed-loop system to follow the reference profile~$x_\text{ref}(a,t)$ by providing a Control Lyapunov Functional. In addition, the proof in the next chapter shows that the overshoot in~\eqref{thm:Linf} as well as the constants~$\alpha_1,\alpha_2$ involved in the CLF~\eqref{CLF_V} are independent of a particular reference trajectory. However, the convergence rate~$L(y_\text{ref})$ in~\eqref{thm:Linf} and \eqref{thm:dotV} is governed by~$y_\text{ref}$.
%
%
%
%
\section{Proof of the Main Results} \label{sec:ProofMain}
%
Before starting with the proof of the main theorem, we compute the closed-loop dynamics in dependence of the delay variables $\delta(t)$ and $\eta(t)$ and $\bm{z}(t)$. 
%
%
We use the control law~\eqref{D_FF}--\eqref{D_CL}, ODE~\eqref{ODE_eta} and definition \eqref{Y_psidelta} and conclude to%
%
%
\begin{align}
\label{eta_dot}
\dot{\eta}(t) 		=&D^*~~\,-\frac{\dot{y}_\text{ref}(t)}{y_\text{ref}(t)}-\text{sat}\Big(z_2(t)-\frac{\dot{y}_\text{ref}(t)}{y_\text{ref}(t)}+\gamma\eta(t)+\gamma\delta(t)\Big)\\
\label{z1_dot}
\dot{z}_1(t)    =&z_2(t)-\frac{\dot{y}_\text{ref}(t)}{y_\text{ref}(t)}-\text{sat}\Big(z_2(t)-\frac{\dot{y}_\text{ref}(t)}{y_\text{ref}(t)}+\gamma\eta(t)+\gamma\delta(t)\Big)\nonumber\\
												&-l_1[z_1(t)-\eta(t)-\delta(t)]\\
\label{z2_dot}
\dot{z}_2(t)   =&-l_2[z_1(t)-\eta(t)-\delta(t)].
%
%
%
\end{align}
Recall that $\psi(t-a)$ and the mapping $\delta(t)$ are input independent by virtue of~\eqref{IC_etapsi} and \eqref{def_delta}. With this result, we are in the position to describe the closed-loop dynamics of the PDE system~\eqref{PDE}--\eqref{OE} completely with the delay model:
%
%
%
\begin{lemma}[controlled dynamics] \label{lem:delay_cl}
Consider the controlled delay model governed by~\eqref{IDE_psi},\,\eqref{eta_dot}--\eqref{z2_dot} with the ICs \eqref{IC_etapsi}, and $\bm{z}_0\in\mathbb{R}^2$. The unique solution~$\eta(t)$, $\psi(t-a)$ and $\bm{z}(t)$ exists for all~$t\geq0$, $a\in[0,A]$ and parameterizes the solution of the PDE model~\eqref{PDE}--\eqref{OE},\,\eqref{D_FF}--\eqref{D_CL} by virtue of \eqref{x_psieta}--\eqref{y_psieta}.
\end{lemma}
\textbf{Proof of Lemma~\ref{lem:delay_cl}.}
Assumptions (A$_1$)--(A$_2$) ensure the existence of a unique solution of~\eqref{PDE}--\eqref{OE} as elucidated in Lemma~3.2 of~\cite{KaraKrstic.2016}. The solution $\psi$ of~\eqref{IDE_psi} belongs to $\mathcal{C}^1(\mathbb{R}^+;\mathbb{R})$, since it coincides with the solution of the delay differential equation\\[-3ex]
\begin{align}
\label{DDE_psi}
\frac{\text{d}\psi(t)}{\text{d}t}=\tilde{k}(0)\psi(t)-\tilde{k}(A)\psi(t-A)+\int_0^A \tilde{k}'(a) \psi(t-a)\text{d}{a}\\[-5ex]\nonumber
\end{align}
with the same IC. On this basis, the function defined by~\eqref{x_psieta} is of class $\mathcal{X}$ for any $x_0\in\mathcal{X}$, $y_\text{ref}\in\mathcal{Y}$ and any piecewise continuous input $D:~\mathbb{R}^+_0\to[D_\text{min},D_\text{max}]$. Define the set
$
\Omega = \{(a,t)\in (0,A)\times\mathbb{R}^+:~(a-t)\notin B\cup\{0,A\} \},
$
where $B\subseteq(0,A)$ is the finite (possibly empty) set where the derivative of $x_0\in\mathcal{X}$ is not defined or is not continuous.\\
The function~$\delta(t)\in\mathbb{R}$ is continuous and well-defined by \eqref{def_delta} for all $\psi(t-a)$ and every kernel $g(a)$ defined by \eqref{g(a)} under assumption (A$_1$).
As a result, the solution of the differential equations \eqref{eta_dot}--\eqref{z2_dot} is unique and exists locally for every $(\eta_0,\bm{z}_0)\in\mathbb{R}\times\mathbb{R}^2$. Since the right-hand sides of \eqref{eta_dot}--\eqref{z2_dot} satisfy a linear growth condition, the solutions $\bm{z}(t)$ and $\eta(t)$ exist for all $t\geq0$.
The next step is to prove that the solutions of \eqref{eta_dot}--\eqref{z2_dot} can be utilized for parameterizing the PDE solution. For all $t\geq0$ and $a\in[0,A]$ the relations
%
%
%
\begin{align}
\label{phi_eta}
\eta(t)=\ln \frac{\Pi\big(\Phi_x[t]\big)}{y_\text{ref}(t)},~~\psi(t-a)=\frac{\Phi_x(a,t)}{x^*(a)\Pi\big(\Phi_x[t]\big)}-1
%
\end{align}
hold true from \eqref{IC_etapsi}--\eqref{x_psieta} with $\Pi(x^*(a)\psi(t-a))=0$. In addition, it is straightforward to show the following identities on $(a,t)\in\Omega$ using~\eqref{ODE_eta}--\eqref{x_psieta}
\begin{align}
\label{phi_PDE}
\frac{\partial \Phi_x(a,t)}{\partial t}+\frac{\partial \Phi_x(a,t)}{\partial a}=-[\mu(a)+D(t)]\Phi_x(a,t).
\end{align}
In accordance to \eqref{BC},~\eqref{IC} we get
\begin{align}
\Phi_x(0,t) &=\langle k,\Phi_x[t]\rangle&&\text{for all}~t\geq0,\\
\Phi_x(a,0) &=x_0(a) &&\text{for all}~a\in[0,A],
\end{align}
which imply the relation $x(a,t)=\Phi_x(a,t)$ in conjunction with \eqref{phi_eta},~\eqref{phi_PDE}. At last, we evaluate the expression $\langle p, \Phi_x[t]\rangle$ with~\eqref{x*(a)},\,\eqref{g(a)} and obtain~\eqref{y_psieta}, which completes the proof.$\hfill\square$\\[2ex]
For proving Theorem~\ref{thm:GES}, we need further lemmas which characterize the dynamics of the observer error and the scalar ODE-state~$\eta(t)$. The proofs are given in Appendix~\ref{AppA}.
%
%
\begin{lemma}[stability of observer dynamics] \label{lem:J(e)}
Consider the observer error~$\bm{e}(t)$ defined by~\eqref{def_e(t)} in context of the closed loop system~\eqref{IDE_psi},\,\eqref{eta_dot}--\eqref{z2_dot} under the assumptions of Theorem~\ref{thm:GES}. Define the positive quadratic form 
\begin{align}
\label{def_J(e)}
J(\bm{e})=e_1^2-p_1e_1e_2+p_2e_2^2=:\bm{e}^{\intercal}\bm{P}\bm{e}>0
\end{align}
Then, there exist positive constants $K_1, K_2,\beta_1, \beta_2>0$ such that the following inequalities holds for all $t\geq0$
\begin{align}
\label{estimate_e2}
K_1\big|\bm{e}(t)\big|^2 &\leq J\big(\bm{e}(t)\big) \leq K_2\big|\bm{e}(t)\big|^2,\\[1ex]
\label{dJ/dt}
\frac{\emph{d} J\big(\bm{e}(t)\big)}{\emph{d}t}&\leq -2 \beta_1 J\big(\bm{e}(t)\big) +\beta_2 |\delta(t)|^2.
\end{align}
\end{lemma}
\begin{lemma}[stability of scalar dynamics] \label{lem:deta2/dt}
Consider the dynamics of the state~$\eta(t)$ governed by~\eqref{IC_etapsi},\,\eqref{eta_dot}--\eqref{z2_dot} under the assumptions of Theorem~\ref{thm:GES}. Define the positive constant
$
\mu_2  = 8\cdot\gamma^{-1} (D_\emph{max}-D_\emph{min})
$
and the functional 
\begin{align}
\label{mu1(yref)}
\mu_1(y_\emph{ref})=&\min(2,\gamma)\cdot\min\bigg\{1,D^*-D_{\min} - \sup_{\tau\geq0} \frac{\dot{y}_\emph{ref}(\tau)}{y_\emph{ref}(\tau)}, \nonumber\\&D_{\max}-D^* +\inf_{\tau\geq0} \frac{\dot{y}_\emph{ref}(\tau)}{y_\emph{ref}(\tau)}\bigg\}
\end{align}
depending on the reference trajectory. Then, for every~$y_\emph{ref}\in\mathcal{Y}$ the following estimate holds with $\mu_1(y_\emph{ref})>0$ for all $t\geq0$
\begin{align}
\label{estimate_eta2_dot}
\frac{\text{d} \eta^2(t)}{\text{d}t}\leq -\mu_1(y_\emph{ref})\frac{\eta^2(t)}{1+\sqrt{\eta^2(t)}}+\mu_2 \bigg|e_2(t)+\gamma\delta(t)\bigg|.
\end{align}
\end{lemma}
At last, define the functionals
\begin{align}
\label{V_zetapsi}
\widetilde{V}(\bm{e},\eta,\psi) =&~\eta^2+\alpha_1 \sqrt{\widetilde{Q}(\bm{e},\psi)}+\alpha_2\widetilde{Q}(\bm{e},\psi)\\
\label{Q_zetapsi}
\widetilde{Q}(\bm{e},\psi) =&~\frac{M}{2}\left(\frac{\left\|\text{e}^{-\sigma a} \psi(-a)\right\|_{\infty}}{1+\min \{ 0,\min_{a\in[0,A]} \psi(-a)\}}\right)^2\nonumber\\
 &+e_1^2-p_1e_1e_2+p_2e_2^2.
\end{align}
with the constants~$M$, $\alpha_1$, $\alpha_2$ to be selected in the following.\\[1ex]
\textbf{Proof of Theorem~\ref{thm:GES}.} At first we analyze the dynamics of the function~$\delta(t)$, see \eqref{def_delta}, governed by the IDE~\eqref{IDE_psi}. To this end, we make use of results on the stability of IDEs in~\cite{Karafyllis.2013,KaraKrstic.2016}. More specifically, we use Corollary~4.6 in~\cite{KaraKrstic.2016} as a main tool and preliminarily discuss necessary assumptions in order to apply it properly to the IDE subsystem governed by~\eqref{IDE_psi} with IC~\eqref{IC_etapsi}. The existence of a constant $\sigma>0$, satisfying\\[-3ex]
\begin{align}
 \int_0^A \e^{\sigma a} \left|\tilde{k}(a)-\lambda \frac{\int_a^A\tilde{k}(s)\text{d}{s}}{\int_0^A s\tilde{k}(s)\text{d}{s}} \right|\text{d}{a}<1
\end{align}
is a direct consequence of assumption~(B$_3$) with a constant $\lambda>0$. Next, consider the continuous functional defined by
\begin{align}
P(\psi_t)=\frac{1}{\int_0^A s\tilde{k}(s)\text{d}{s}}\int_0^A \int_a^A\psi_t(-a)\tilde{k}(s)\text{d}{s}\text{d}{a}
\end{align}
using the notation $\psi_t(-a)=\psi(t-a)$ and notice the fact $P(\psi_0)=0$  by virtue of the IC~\eqref{IC_etapsi}. Therefore, it follows from Corollary~4.6 in~\cite{KaraKrstic.2016} that
\begin{align}
\label{W(psi)}
W(\psi_t) = \left\|\e^{-\sigma a} \psi_t(-a)\right\|_{\infty} 
\end{align}
is a Lyapunov functional of the delay subsystem~\eqref{IC_etapsi},~\eqref{IDE_psi} by means of the differential inequality $\dot{W}^+(\psi_t)\leq - \sigma W(\psi_t)$ for all $t\geq0$ with $\sigma>0$. In addition, Remark 4.7~\cite{KaraKrstic.2016} implies that the functional defined by
\begin{align}
\label{C(psi)}
C(\psi_t) = 1+\min \left( 0,\min_{a\in[0,A]} \psi_t(-a)\right)
\end{align}
is non-increasing. As a result, it follows $\psi_t(-a)>-1$ for all $a\in[0,A]$ and $t\geq0$. Using the fact $|\ln(1+f)|\leq|f|/ (1+\min(f,0))$ for all $f>-1$, we conclude from~\eqref{def_delta}:
\begin{equation}
\label{ineq_delta}
|\delta(t)|\leq \frac{\text{e}^{\sigma A}\left\|\text{e}^{-a\sigma} \psi_t(-a)\right\|_{\infty}}{1+\min\{0,\min_{a\in[0,A]} \psi_t(-a)\}}~\text{for all }t\geq0.
\end{equation}
%
The second part of the current proof aims at deriving bounds for the analysis parameters~$\alpha_1$ and~$\alpha_2$ independent of a reference trajectory, such that~$V(\bm{z},x,t)$ is a CLF of the closed loop~\eqref{PDE}--\eqref{OE}, \eqref{D_FF}--\eqref{D_CL}. For this purpose, the following relations between the delay model~\eqref{IC_etapsi}--\eqref{y_psieta} and the PDE problem~\eqref{PDE}--\eqref{OE} are a direct consequence of definition~\eqref{V_zetapsi}, \eqref{Q_zetapsi}:
\begin{align}
\label{Vtilde=V}
\widetilde{V}\big(\bm{e}(t),\eta(t),\psi_t\big) &=V\big(\bm{e}(t)+[\eta(t),D^*]^\intercal,\Phi_x(a,t),t\big) \\
\widetilde{Q}\big(\bm{e}(t),\psi_t\big)&=Q\big(\bm{e}(t)+[\eta(t),D^*]^\intercal,\Phi_x(a,t),t\big)
\end{align}
Subsequently, we obtain the following inequality for the Dini-derivative of~\eqref{Q_zetapsi} for any~$M>0$
\begin{align}
\dot{\widetilde{Q}}\,^+\big(\bm{e}(t),\psi_t\big)&\leq - 2 \beta_1 J\big(\bm{e}(t)\big)-\left(\sigma M - \beta_2\text{e}^{2\sigma A}\right)\nonumber\\ 
&~~\times \left(\frac{\left\|\text{e}^{-\sigma a} \psi_t(a)\right\|_{\infty}}{1+\min \{ 0,\min_{a\in[0,A]} \psi_t(a)\}}\right)^2.
\end{align}
by identifying the non-increasing functional~\eqref{C(psi)} in the denominator of~\eqref{Q_zetapsi}.
Choose~$M$ such that $M\sigma>\beta_2\text{e}^{2\sigma A}$ holds. To this end, we conclude that the estimate
\begin{equation}
\label{estimate_Qdot}
\dot{\widetilde{Q}}\,^+\big(\bm{e}(t),\psi_t\big)\leq -2 \beta  \widetilde{Q}\big(\bm{e}(t),\psi_t\big)
\end{equation}
holds with the constant $\beta=\min\left(\beta_1,\sigma-\text{e}^{2\sigma A}\beta_2M^{-1}\right)>0.$
%
%
%
%
Beyond that, for all $t\geq0$ with $\widetilde{Q}(\bm{e}(t),\psi_t) > 0$ we combine
\begin{align}
&\frac{\widetilde{Q}^{\nicefrac{1}{2}}(\bm{e},\psi)\big|_{t+h}-\widetilde{Q}^{\nicefrac{1}{2}}(\bm{e},\psi)\big|_{t}}{h} \leq\frac{\widetilde{Q}(\bm{e},\psi)\big|_{t+h} -\widetilde{Q}(\bm{e},\psi)\big|_{t}}{h} \nonumber\\&~~~~~~~~~~~~~~~~~~~~~\times{\left(\widetilde{Q}^{\nicefrac{1}{2}}(\bm{e},\psi)\big|_{t+h}+\widetilde{Q}^{\nicefrac{1}{2}}(\bm{e},\psi)\big|_{t}\right)^{-1}}
\end{align}
with the fact $\lim_{h\to0}  \widetilde{Q}^{\nicefrac{1}{2}}({e},\psi)\big|_{t+h} = \lim_{h\to0}  \widetilde{Q}^{\nicefrac{1}{2}}({e},\psi)\big|_{t}$ and get the differential inequality
\begin{equation}
\label{estimate_sqrtQdot}
\dot{\sqrt{\widetilde{Q}\big(\bm{e}(t),\psi_t\big)}}^+\leq - \beta  \sqrt{\widetilde{Q}\big(\bm{e}(t),\psi_t\big)}.
\end{equation} 
Moreover, \eqref{estimate_Qdot} implies that the function $t\mapsto\widetilde{Q}(\bm{e}(t),\psi_t)$ is non-increasing, so if we have $\widetilde{Q}(\bm{e},\psi)|_t=0$ we also have $\widetilde{Q}^{\nicefrac{1}{2}}(\bm{e},\psi)|_{t+h}=0$ for every $h\geq0$, because $ \widetilde{Q}(\bm{e},\psi)\big|_{t+h}\leq\widetilde{Q}(\bm{e},\psi)\big|_{t}=0$. We conclude that \eqref{estimate_sqrtQdot} holds for all $t\geq0$ with $\widetilde{Q}(\bm{e}(t),\psi_t)=0$. Thus, \eqref{estimate_sqrtQdot} is valid for all $t\geq0$.\\[1ex]
With the results up to this point, we are in the position to establish the following inequality for the Dini-derivative of~\eqref{V_zetapsi} for all $t\geq0$ based on \eqref{dJ/dt}, \eqref{estimate_eta2_dot}, \eqref{estimate_Qdot}, \eqref{estimate_sqrtQdot}, viz.
\begin{align}
\dot{\widetilde{V}}\,^+\big(\bm{e}(t),\eta(t),\psi_t\big) \leq& -\mu_1(y_\text{ref})\frac{\eta^2(t)}{1+\sqrt{\eta^2(t)}}+\mu_2\big|e_2(t)\big|\nonumber\\&+\mu_2\gamma\big|\delta(t)\big|- \alpha_1 \beta  \sqrt{\widetilde{Q}\big(\bm{e}(t),\psi_t\big)}\nonumber\\[1ex]& - 2\alpha_2\beta\widetilde{Q}\big(\bm{e}(t),\psi_t\big).
\end{align}
Next, we make use of the inequalities  
\begin{align}
\sqrt{K_1}|e_2(t) |  \leq \sqrt{K_1} |\bm{e}(t) |&\leq \sqrt{\widetilde{Q}\big(\bm{e}(t),\psi_t\big)}\\
\frac{\sqrt M }{\sqrt{2}\text{e}^{\sigma A}}|\delta(t)|&\leq \sqrt{\widetilde{Q}\big(\bm{e}(t),\psi_t\big)},
\end{align}
which are a consequence of  \eqref{estimate_e2}, \eqref{ineq_delta} and definition~\eqref{Q_zetapsi}, in conjunction with inequality~\eqref{estimate_Qdot} ending up with
\begin{align}
\label{eq52}
\dot{\widetilde{V}}\,^+\big(\bm{e}(t),\eta(t),\psi_t\big) \leq & -\frac{\mu_1(y_\text{ref})\eta^2(t)}{1+\sqrt{\eta^2(t)}}- 2\alpha_2\beta\widetilde{Q}\big(\bm{e}(t),\psi_t\big) \nonumber\\&-\bigg(\alpha_1 \beta
-\frac{8\sqrt{2}(D_\text{max}-D_\text{min})\text{e}^{\sigma A}}{{\sqrt M}} \nonumber\\
& - \frac{8(D_\text{max}-D_\text{min})}{\gamma \sqrt{K_1}}\bigg) \sqrt{\widetilde{Q}\big(\bm{e}(t),\psi_t\big)}.
\end{align}
We select the analysis parameters $\alpha_1$ and $\alpha_2$ independent of the reference trajectory $y_\text{ref}\in\mathcal{Y}$ as
\begin{align}
\alpha_2 > 0,~~~
\alpha_1 > \frac{8(D_\text{max}-D_\text{min})}{\beta}\left[\frac{1}{\gamma \sqrt{K_1}}+\frac{\sqrt{2}\text{e}^{\sigma A}}{{\sqrt M}}\right]
\end{align}
so that the differential inequality 
\begin{align}
\label{proof:dotV}
\dot{\widetilde{V}}\,^+\big(\bm{e}(t),\eta(t),\psi_t\big) &\leq - L(y_\text{ref}) \frac{\widetilde{V}\big(\bm{e}(t),\eta(t),\psi_t\big)}{1+ \sqrt{ \widetilde{V}\big(\bm{e}(t),\eta(t),\psi_t\big)}}
\end{align}
holds for all $t\geq0$ and $(\bm{e}(t),\eta(t),\psi_t(-a))\in\mathbb{R}^4\setminus\{\textbf{0}\}$ with
\begin{align}
\label{def_L(yref)}
L(y_\text{ref})=&~\min\bigg(\beta-  \frac{8(D_\text{max}-D_\text{min})}{\alpha_1} \left[  \frac{1}{\gamma\sqrt{K_1}}+\frac{\sqrt{2}\text{e}^{\sigma A}}{{\sqrt M}}\right],\nonumber\\
&~\min(2,\gamma)\cdot\min\bigg\{1,D^*-D_{\min} - \inf_{\tau\geq0} \frac{\dot{y}_\text{ref}(\tau)}{y_\text{ref}(\tau)}, D_{\max}\nonumber\\&~-D^* +\sup_{\tau\geq0} \frac{\dot{y}_\text{ref}(\tau)}{y_\text{ref}(\tau)}\bigg\}\bigg) >0.
\end{align}
Finally, \eqref{thm:dotV} is a consequence of~\eqref{mu1(yref)},\,\eqref{Vtilde=V},\,\eqref{eq52}--\eqref{proof:dotV}. In addition, \eqref{def_L(yref)} shows how the control gain~$\gamma$, the reference trajectory~$y_\text{ref}$ and the input constraint $D(t)\in[D_{\min},D_{\max}]$ influence the convergence rate.\\
According to Lemma 5.2 in~\cite{KaraKrstic.2016} the following estimate bounds the solution of the differential inequalities~\eqref{thm:dotV}, \eqref{proof:dotV}
\begin{align}
\label{V_GES}
\widetilde{V}\big(\bm{e}(t),\eta(t),\psi_t\big)  \leq \text{e}^{-\frac{L(y_\text{ref})}{2}t} \cdot\widetilde{V}_0 \text{e}^{\max(0,\widetilde{V}_0 -1)} 
\end{align}
with the abbreviation $\widetilde{V}_0=\widetilde{V}(\bm{e}_0,\eta_0,\psi_0)$ and is employed to prove~\ref{thm:GES} in the following. In order to verify the existence of a function~$\kappa\in\mathcal{K}_\infty$ take the following estimate as a basis which again uses the inequality $|\ln(1+f)|\leq |f| /(1+\min(0,f))$ for all $f>-1$ and recall the relation~$x(a,t)=\Phi_x(a,t)$ by virtue of~\eqref{x_psieta} and \eqref{ineq_delta}.
\begin{align}
\left|\ln\frac{x(a,t)}{x_\text{ref}(a,t) }\right|\leq |\eta(t)|+\frac{\text{e}^{\sigma A}\left\|\text{e}^{-a\sigma} \psi_t(a)\right\|_{\infty}}{1+\min\{0,\min_{a\in[0,A]} \psi_t(a)\}}
\end{align}
Since the statement holds for all $a\in[0,A]$ it is also fulfilled for the $L^\infty$-norm of the left-hand side. With this fact, \eqref{V_GES} and the inequalities for $t\geq0$
\begin{align}
\label{estimate_components}
 \eta^2(t) &\leq \widetilde{V}\big(\bm{e}(t),\eta(t),\psi_t\big) \\ \alpha_1\sqrt{\widetilde{Q}\big(\bm{e}(t),\psi_t)}
 &\leq \widetilde{V}\big(\bm{e}(t),\eta(t),\psi_t\big) ,
\end{align}
which result from~\eqref{V_zetapsi}, we get for all $t\geq0$:
\begin{align}
\label{Linf_estimate_V_0}
%
%
%
%
\underbrace{\left\|\ln\frac{x(a,t)}{x_\text{ref}(a,t)}\right\|_{\infty}}_{=:\varsigma(t)}+|\bm{e}(t)|\leq& \Bigg(\widetilde{V}_0^{\nicefrac{1}{2}}+\frac{\sqrt{2}\e^{\sigma A}}{\alpha_1\sqrt{M}}\widetilde{V}_0+\frac{\widetilde{V}_0}{\sqrt{K_1}\alpha_1}\Bigg)\nonumber \\[-3ex]
&\times\text{e}^{\max(0,\widetilde{V}_0 -1)}\text{e}^{-\frac{L(y_\text{ref})}{4} t} 
\end{align} 
If $M,\alpha_1>0$ are sufficiently large constants satisfying~$\alpha_1\min(\sqrt{K_1},\sqrt{M/2})\geq2$, it follows with $\widetilde{\kappa}\in\mathcal{K}_\infty$
\begin{align}
\varsigma(t)+|\bm{e}(t)|&\leq \text{e}^{-\frac{L(y_\text{ref})}{4} t} \cdot \text{e}^{\sigma A} \left(\widetilde{V}_0^{\nicefrac{1}{2}}+\widetilde{V}_0\right)\text{e}^{\max(0,\widetilde{V}_0 -1)} \\   
												&=\text{e}^{-\frac{L(y_\text{ref})}{4} t} \cdot \widetilde{\kappa}\left(\widetilde{V}_0\right)
\end{align}
Next, we show the existence of a function~$\kappa_V\in\mathcal{K}_\infty$ (independent of $y_\text{ref}$) such that $\widetilde{V}_0\leq \kappa_V(\varsigma_0+|\bm{e}_0|)$ where $\varsigma_0=\varsigma(0)$. For this purpose, it suffices to show that there exist two functions $\kappa_\eta,  \kappa_\psi\in\mathcal{K}_\infty$ which bound the ICs~\eqref{IC_etapsi}
%
%
%
%
\begin{align}
\label{kappa_eta}
|\eta_0| &\leq \kappa_\eta(\varsigma_0)\\
\label{kappa_psi}
\frac{\left\|\text{e}^{-a\sigma} \psi_0(a)\right\|_{\infty}}{1+\min\{0,\min_{a\in[0,A]} \psi_0(a)\}}&\leq \kappa_\psi(\varsigma_0).
\end{align}
Indeed, if \eqref{kappa_eta}--\eqref{kappa_psi} hold then definition~\eqref{estimate_e2} and \eqref{Q_zetapsi} imply 
\begin{align}
\label{kappa_Q}
\widetilde{Q}(\bm{e}_0,\psi_0)  \leq K_2 |\bm{e}_0|^2+\frac{M}{2}\big(\kappa_\psi(\varsigma_0)\big)^2\leq \kappa_Q(\varsigma_0+|\bm{e}_0|)
\end{align}
where $\kappa_Q(z)=(K_2+M/2)(z+\kappa_\psi(z))^2$ for any argument $z\geq0$. We get the desired function $\widetilde{V}_0\leq \kappa_V(\varsigma_0+|\bm{e}_0|)$ with definition~\eqref{V_zetapsi} and inequalities~\eqref{kappa_eta}--\eqref{kappa_Q}, where
\begin{align}
\label{kappa_V}
\kappa_V(z)=\kappa^2_\eta(z)+\alpha_1 \sqrt{\kappa _Q(z)} +\alpha_2 \kappa_Q(z),~~z>0.
\end{align}
Furthermore, we identify $|\eta_0|\leq\kappa_\eta(\varsigma_0)=\varsigma_0$ by virtue of the definitions~\eqref{def_Pi},\,\eqref{IC_etapsi} for all $x_0\in\mathcal{X}$
\begin{align}
\text{e}^{-\varsigma_0}\leq \min_{a\in[0,A]} \frac{x_0(a)}{x_\text{ref}(a,0)} \leq\frac{\Pi(x_0)}{y_\text{ref}(0)}\leq \max_{a\in[0,A]} \frac{x_0(a)}{x_\text{ref}(a,0)}\leq\text{e}^{\varsigma_0}.
\end{align}
The following inequalities for~$\psi_0(a)$ hold in same manner
\begin{align}
\label{estimate_psi0_varsigma0}
\text{e}^{-2\varsigma_0}\leq\psi_0(a)+1 \leq|\psi_0(a)|+1\leq \text{e}^{2\varsigma_0}.
\end{align}
We get the result $\kappa_\psi(\varsigma_0)=\text{e}^{2\varsigma_0}(\text{e}^{2\varsigma_0}-1)$ as a valid bound for all $\varsigma_0\geq0$ when plugging \eqref{estimate_psi0_varsigma0} into~\eqref{kappa_psi}. With the estimate
%
\begin{align}
\varsigma(t)+|\bm{e}(t)|\leq \widetilde{\kappa}\big(\kappa_V(\varsigma_0+|\bm{e}_0|)\big)\cdot \text{e}^{-\frac{L(y_\text{ref})}{4} t}
\end{align}
we finally conclude that \eqref{thm:Linf} holds for $t\geq0$ with $\kappa=\widetilde{\kappa}\circ\kappa_V$. The proof is complete.$\hfill\square$
%
%
%
%
\section{Illustrative Example} \label{sec:Simulation}
%
Applying the method of weighted residuals and using the age-weight $w(a)$~we convert PDE~\eqref{PDE} to a weak form
\begin{align}
\label{PDE_weak}
\int_0^A \bigg[\frac{\partial x(a,t)}{\partial t}+&\frac{\partial x(a,t)}{\partial a}\bigg]w(a)\text{d}{a} =\nonumber\\
&-\int_0^A\Big[[\mu(a)+D(t)]x(a,t)\Big]w(a)\text{d}{a} 
\end{align} 
with the goal of reducing the IBVP~\eqref{PDE}--\eqref{OE} to an approximate finite-dimensional initial-value problem for simulation purposes. We therefore use the ansatz $\tilde{x}(a,t)=\bm{\varphi}^{\intercal}(a)\bm{\lambda}(t)$ with a set of $N$ BC-compatible and linear independent trial functions $\varphi_k\in P\mathcal{C}^1([0,A])$ and the temporal weights~$\bm{\lambda}(t)\in\mathbb{R}^{N}$ are introduced. In order to minimize the error with respect to the exact solution~$x(a,t)$, we choose the weights to consist of a linear combination of the trial functions $w(a)=\bm{w}^{\intercal}\bm{\varphi}(a)$ in Galerkin-manner. Taking this procedure as a basis, we obtain the following residual equation from~\eqref{PDE_weak} 
\begin{align}
0=~&\bm{w}^{\intercal}\int_0^A \bm{\varphi}(a)\bigg[\bm{\varphi}^{\intercal}(a)\dot{\bm{\lambda}}(t)+\left(\bm{\varphi}'\right)^{\intercal}(a)\bm{\lambda}(t)\nonumber\\ &+[\mu(a)+D(t)]\bm{\varphi}^{\intercal}(a)\bm{\lambda}(t)\bigg]\text{d}{a}.
\end{align}
For all $\bm{w}\neq\textbf{0}$ it is possible to extract an initial-value problem for $t\geq0$ with the matrices~$\bm{M},\bm{N}\in\mathbb{R}^{N\times N}$
\begin{align}
\bm{M}&=~~~\int_0^A \bm{\varphi}(a)\bm{\varphi}^{\intercal}(a)\text{d}{a}\\
\bm{N}&=-\int_0^A \bm{\varphi}(a)\left(\bm{\varphi}'\right)^{\intercal}(a)+\mu(a)\bm{\varphi}(a)\bm{\varphi}^{\intercal}(a)\text{d}{a}
\end{align}
consisting of the ODE
\begin{align}
\label{lambda_dot}
\dot{\bm{\lambda}}(t)=\left[\bm{M}^{-1}\bm{N}-\bm{I} D(t)\right]\bm{\lambda}(t),~~~t>0
\end{align}
and the IC~$\bm{\lambda}(0)=\bm{\lambda}_0=[1~0~\cdots~0]^{\intercal}$ for the specific choice~$\varphi_1=x_0\in\mathcal{X}$, if the initial profile is linear independent of the remaining trial functions. These are adopted similar to a Karhunen-Lo{\`e}ve-decomposition. Firstly, $\varphi_2(a)=x^*(a)$ guarantees accurate asymptotic properties, since the steady-state solution of \eqref{PDE}--\eqref{OE} is included in~$\bm{\varphi}(a)$. For $k\geq2$ we choose with $j\in[2,\frac{N}{2}]$
\begin{align}
\varphi_{2j-1}(a)           &=&\cos(\omega_k a) \text{e}^{\sigma_k a} \varphi_2(a)&\qquad\qquad~~~~~~~~\\
\varphi_{2j}(a)      &=&\sin(\omega_k a)  \text{e}^{\sigma_k a}\varphi_2(a)&.
\end{align}
with the remaining complex conjugated solutions $\sigma_k\pm j\omega_k$ of $\int_0^A  \tilde{k}(a) \text{e}^{ \sigma_ka\pm j\omega_ka}\text{d}{a}=1$. Note that $\bm{M}>0$ is a direct consequence of the linear independent trial functions. We checked the solution such that the positivity condition~$\hat{x}(a,t)>0$ is not violated at any time and age. At last, the output of the finite-order approximation is given by
\begin{figure}
\centering
\includegraphics[width=0.48\textwidth]{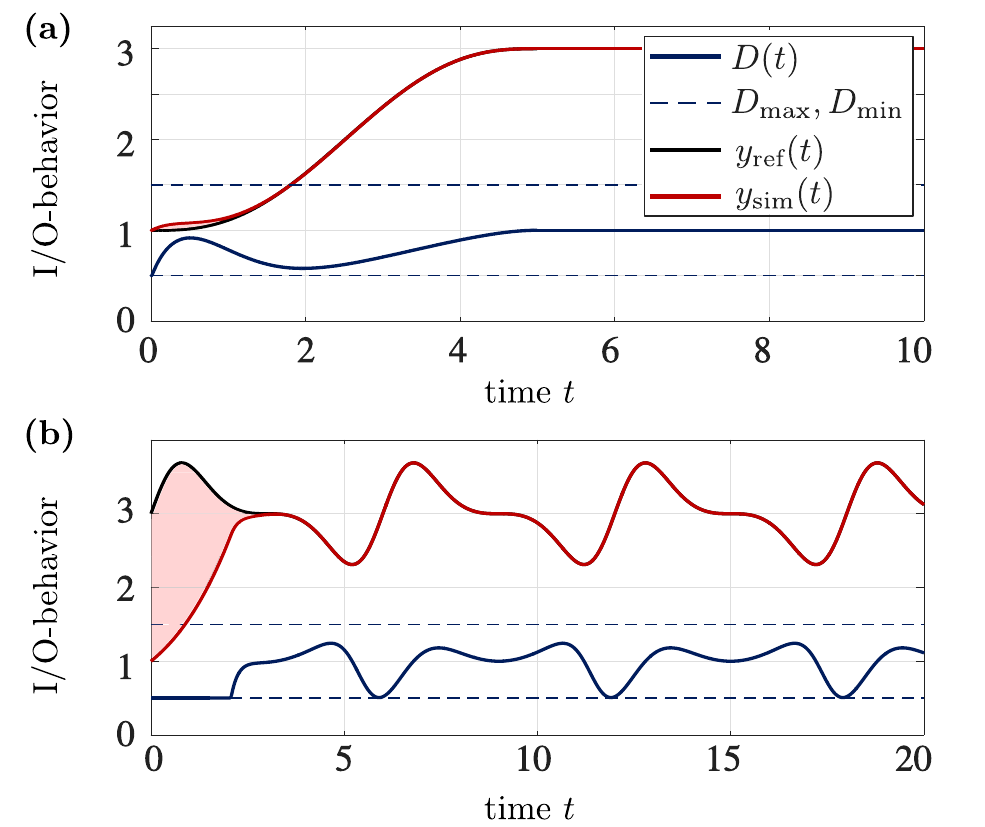}  
\caption{Closed loop I/O-behavior of a set point change (a) governed by~\eqref{yref_tran} and the periodical trajectory~\eqref{yref_peri} in plot~(b).}   
\label{fig:IO_simulation}                            
\end{figure}
\begin{align}
{y}_\text{sim}(t)=\int_0^A p(a) \bm{\varphi}^{\intercal}(a)\text{d}{a}\bm{\lambda}(t) = \bm{p}^{\intercal} \bm{\lambda}(t)
\end{align}
\textit{Trial System. } In order to demonstrate the introduced asymptotic tracking controller, consider the following trial system in the age-domain up to $A=2$ with a symmetric input constraint $D(t)\in[D_\text{min},D_\text{max}]=[0.5,1.5]$ around the equilibrium value $D^*=1$. We use the mortality rate~$\mu(a)=\mu=0.1$ and a quadratic motherhood kernel
$
k(a)=k_0 a\cdot (A-a)
$
satisfying the Lotka-Sharpe condition~\eqref{LotkaSharpe}, since
\begin{align}
k_0 &= \left(\int_0^A a (A-a) \text{e}^{-(D^*+\mu)a} \text{d}{a}\right)^{-1}=2.00.
\end{align}
Moreover assume that the yield production rate is proportional to the overall population, namely $p(a)=1$ and the initial profile~$x_0(a)=\varphi_1(a)=-0.054a+\text{e}^{-1.30a}$, by evaluating $\langle1,x_0\rangle=1$. For realizing~$N=6$ we obtain the conjugated pairs $\sigma_k\pm j \omega_k=\{-2.02\pm 4.41j,-2.50,\pm7.62j\}$ computed with a Newton-Raphson scheme.\\[1ex]
Next, the following set of desired trajectories is considered in order to demonstrate control structure \eqref{D_FF}--\eqref{D_CL} with $\gamma=2$, $\bm{l}=[4~8]^{\intercal}$ and the initial guess $D^*_0=D_\text{min}=0.5$
\begin{figure}
\centering
\includegraphics[width=0.48\textwidth]{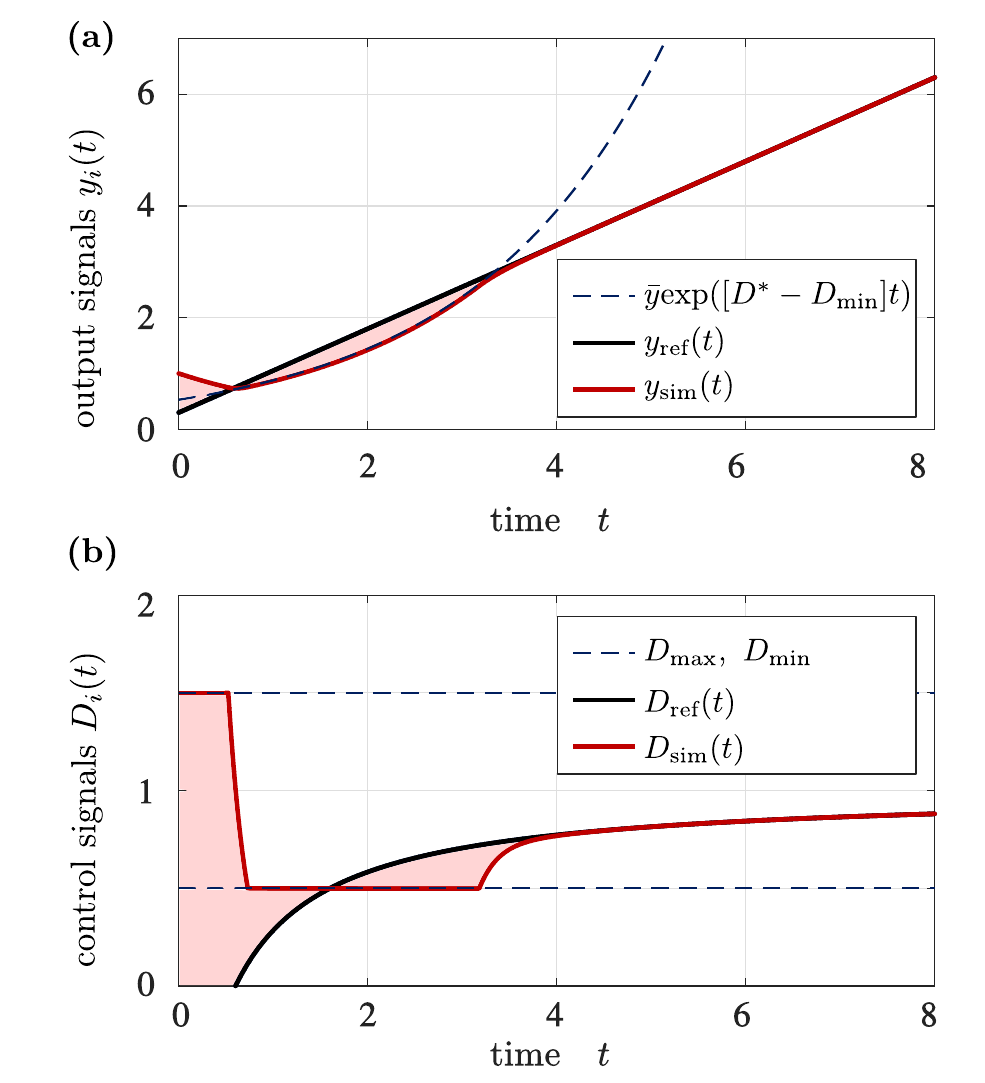}  
\caption{Simulation of the closed loop for a ramp-like reference trajectory~\eqref{yref_ramp} with $y_\text{ref,1}\notin\mathcal{Y}$.}
\label{fig:IO_Yramp}                            
\end{figure}
\begin{align}
\label{yref_ramp}
y_\text{ref,1}(t)&=y_4+y_1 t\\
\label{yref_peri}
y_\text{ref,2}(t) &=y_2 + y_3\sin[\omega t + \sin(\omega t)]\\
\label{yref_tran}
y_\text{ref,3}(t) &=\begin{cases} 
y_0+(y_\Delta-y_0)\sum_{i=1}^3 g_i \left(\frac{t}{t_\Delta}\right)^{i+2}   &,~t\in[0,t_\Delta]\\
y_\Delta 																																																	&,~t>t_\Delta.
 \end{cases}
\end{align}
The parameters are chosen as $y_0=y(0)=1$, $y_\Delta=3$, $y_1=0.75$, $y_2=0.79$, $y_3=0.625$, $y_4 = 0.3$ and $\omega=\frac{2\pi}{6}$. In addition it is straightforward to verify that~$g_i =\{10,~-15,~6 \}$ ensures the desired smoothness of $y_\text{ref,3}\in\mathcal{C}^2$.\\ 
The simulation results are given in Figure~\ref{fig:IO_simulation} showing the I/O-behavior of the closed-loop for the reference trajectories $y_\text{ref,2/3}\in\mathcal{Y}$. 
With \eqref{PDE} and \eqref{lambda_dot} we evaluate the accuracy of the numerical approach on the basis of the residual 
\begin{align}
R(a,t) =&(\bm{\varphi}')^{\intercal}(a)\bm{\lambda}(t)+\bm{\varphi}^{\intercal}(a)\left[\bm{M}^{-1}\bm{N}-\bm{I} D(t)\right]\bm{\lambda}(t)\nonumber\\&+[\mu(a)+D(t)]\bm{\varphi}^{\intercal}(a)\bm{\lambda}(t)
\end{align}
and its $L^2$-norm $r^2(t)=\langle R[t],R[t]\rangle$.
When simulating the results for~$y_\text{ref,3}(t)$, the numerical approach has a mean error~$\bar{r}=4.9\%$ and an steady error~$r(t_\Delta)=1.9\cdot 10^{-9}$ for $t\geq t_\Delta$.\\
Moreover, Figure~\ref{fig:IO_Yramp} shows the simulation results for the ramp-like trajectory~$y_\text{ref,1}\notin\mathcal{Y}$ which is not valid for all $t\geq0$ since the supremum 
\begin{align}
\sup_{t\geq0} \frac{\dot{y}_\text{ref,1}(t)}{y_\text{ref,1}(t)} = \sup_{t\geq0} \frac{y_1}{y_4+y_1t} = 2.5 > D^*-D_\text{min} = 0.5
\end{align}
violates Assumption~(B$_2$). However, the reference trajectory enters the valid area at $t=t_\text{crit}=1.6$ and the tracking with an exponential decay by virtue of Theorem~\ref{thm:GES} holds for $t>t_\text{crit}$. In addition, the situation with an active saturation is emphasized by plot~(b) in Figure~\ref{fig:IO_Yramp} for $t\leq t_\text{crit}$. As long as the control resides at the lower bound $D(t)=D_\text{min}$, the output~$y(t)$ is proportional  to the limiting exponential function 
$
y_\text{min}(t)=\text{e}^{(D^*-D_\text{min})t},
$
and vice versa for $D(t)=D_\text{max}$.
%
\section{Concluding Remarks} \label{sec:Conclusion}
Control of age-structured populations derives its motivation from mathematical demography, pharmaceutical industry and other fields. Using a PDE model for age-structured chemostats makes it possible to capture transient dynamics more faithfully than with ODEs, but demands more care in control design.\\[1ex]
We applied Lyapunov-methods and recent results on IDEs for establishing a global asymptotic tracking with respect to desired yield trajectories and to provide CLFs for the closed loop system. Furthermore, our control design treats state constraints and input saturation in an explicit way and is robust against parameter uncertainties, not requiring the exact equilibrium dilution rate. In our simulation we demonstrate the attractivity even for the worst initial guess on the boundary of the dilution rate's valid interval. With a 2-DOF control structure we decoupled the feedforward and the feedback control. This configuration yields to an enhanced tracking behavior, while the feedback controller can be used to assign the desired stability properties. Control design for an advanced model that includes the substrate dynamics is still an open topic for further research.\\[1ex]
\textit{Acknowledgements.} The authors would like to thank the anonymous reviewers who read the initial version of the paper carefully.

\appendix
\section{Analysis of Closed-Loop Dynamics}\label{AppA}   
\textbf{Proof of Lemma~\ref{lem:J(e)}.}
Consider the dynamic system consisting of a static equation for equilibrium dilution rate $D^*$ and~\eqref{ODE_eta}
\begin{align}
\label{obssys_1}
\dot{\eta}(t)= D^*- \frac{\dot{y}_\text{ref}(t)}{y_\text{ref}(t)}-D(t), \quad
\dot{D}^* = 0.
\end{align}
It is straightforward to verify that \eqref{z_dot} is an observer for the dynamic system~\eqref{obssys_1} using the logarithmic output~\eqref{Y_psidelta} in order to estimate the state $\eta(t)$ and the equilibrium dilution rate $D^*$. As a consequence of \eqref{z_dot}, \eqref{Y_psidelta}, \eqref{obssys_1} and definition~\eqref{def_delta}, the evolution of the error~\eqref{def_e(t)} is governed by
\begin{equation}
\label{e(t)_dot}
\dot{\bm{e}}(t)=\begin{bmatrix} -l_1& 1 \\-l_2 & 0\end{bmatrix}\bm{e}(t)+\begin{bmatrix} l_1\\l_2\end{bmatrix}\delta(t)=\bm{L}\bm{e}(t)+\bm{l}\delta(t).
\end{equation}
with the Hurwitz matrix $\bm{L}$. Next, take~\eqref{e(t)_dot}, the quadratic form~$J(\bm{e})$ given by~\eqref{def_J(e)} and its derivative
\begin{align}
\label{dJdt_1}
\frac{\text{d} J\big(\bm{e}(t)\big)}{\text{d}t}=-\bm{e}^{\intercal}(t)\widetilde{\bm{P}}\bm{e}(t)+\delta(t) \bm{l}^{\intercal}(\bm{P}^{\intercal}+\bm{P})\bm{e}(t)
\end{align}
into account. The symmetric matrices
\begin{align}
\bm{P}=\begin{bmatrix} 1  & -\frac{p_1}{2} \\ -\frac{p_1}{2} & p_2\end{bmatrix},~
\widetilde{\bm{P}}=\begin{bmatrix} 2l_1-l_2p_1  & l_2p_2-\frac{l_1p_1}{2}-1\\           l_2p_2-\frac{l_1p_1}{2}-1 & p_1 \end{bmatrix},
\end{align}
are positive definite due to inequalities~\eqref{ineq_p12}. Hence, the estimates
$
K_1|\bm{e}|^2\leq \bm{e}^{\intercal}\bm{P}\bm{e}\leq K_2|\bm{e}|^2$ and $
\widetilde{K}_1|\bm{e}|^2\leq \bm{e}^{\intercal}\widetilde{\bm{P}}\bm{e}\leq \widetilde{K}_2|\bm{e}|^2
$
hold with positive constants $K_1,\widetilde{K}_1,K_2,\widetilde{K}_2>0$. With Young's inequality we get the estimate for every $\widetilde{K}_1>0$
\begin{align}
\delta(t) \bm{l}^{\intercal}(\bm{P}^{\intercal}+\bm{P})\bm{e}(t) \leq \frac{|\bm{l}^{\intercal}(\bm{P}^{\intercal}+\bm{P})|^2}{2\widetilde{K}_1}|\delta(t)|^2+\frac{\widetilde{K}_1}{2} |\bm{e}(t)|^2
\end{align}
for the last term of~\eqref{dJdt_1} and conclude that the inequality
\begin{align}
\frac{\text{d} J\big(\bm{e}(t)\big)}{\text{d}t}
&\leq-2\cdot \frac{\widetilde{K}_1}{4K_2} \bm{e}^{\intercal}(t)\bm{P}\bm{e}(t)+\frac{2|\bm{P}\bm{l}|^2}{\widetilde{K}_1} |\delta(t)|^2
\end{align}
holds for all $t\geq0$. It is straightforward to use the previous equation and show that inequality \eqref{dJ/dt} holds with
\begin{align}
\beta_1 &=\frac{\widetilde{K}_1^{\ }}{4{K}_2},~~ 
\beta_2 =\frac{(2l_1-l_2p_1)^2+(l_1p_1-2l_2p_2)^2}{2\widetilde{K}_1}
\end{align}
The proof is complete.$\hfill\square$\\[2ex]
%
\textbf{Proof of Lemma~\ref{lem:deta2/dt}.}
For the subsequent proof we will need the following fact.\\[1ex]
\textbf{Fact}. \textit{For any positive constants $a,b>0$ the inequality}
\begin{align}
\label{ineq_xsatx}
z \cdot\text{sat}_{[-a,b]} (z) \geq \min\left(1,a,b\right) \cdot \frac{ z^2}{1+|z|}.
\end{align}
\textit{holds for all} $z\in\mathbb{R}$.\\[1ex]
\textbf{Proof of Fact.} Distinguish between the cases $z<-a$, {$-a\leq z\leq b$, and $z >b$}. Firstly, we have $|z| a>  \min\left(1,a,b\right)(1+|z|)^{-1}$ , which is a true statement since $a\geq\min\left(1,a,b\right)$ and $|z| \geq z^2(1+|z|)^{-1}$. Secondly, for $z\in[-a,b]$ \eqref{ineq_xsatx} simplifies to  
\begin{equation}
z^2\geq z^2 \cdot \frac{ \min\left(1,a,b\right)}{1+|z|}.
\end{equation}
Due to the fact $1\geq \min\left(1,a,b\right)(1+|z|)^{-1}$ the statement also holds in the present case. Finally, for $z >b$ we receive $z \cdot\text{sat}_{[-a,b]} (z)=z\cdot b$ leading to a true statement for the same reasons like the first case. The proof is complete.$\hfill\square$\\[1ex]
Now, we are ready to prove Lemma~\ref{lem:deta2/dt}. Consider the closed-loop dynamics of the state~$\eta(t)$ governed by~\eqref{eta_dot}. 
%
%
Assumption~(B$_2$) implies $D^*-\dot{y}_\text{ref}(t) y^{-1}_\text{ref}(t)\in(D_\text{min},D_\text{max})$ for all $t\geq0$ as a consequence of~\eqref{ValidYref}. It then follows from~\eqref{eta_dot} with \eqref{def_e(t)} for all $t\geq0$
\begin{align}
\label{dot_eta(t)}
\dot{\eta}(t) &= - \text{sat}_{\mathcal{D}(t)} \Big(e_2(t)+\gamma\eta(t)+\gamma\delta(t) \Big)
\end{align}
with the time-varying interval 
$
\mathcal{D}(t) = [-\bar{D}_{\min}(t) ,\bar{D}_{\max}(t)]
$, where
\begin{align}
\label{def_Dbarmin}
\bar{D}_{\min}(t) &= D^*-D_{\min} - \frac{\dot{y}_\text{ref}(t)}{y_\text{ref}(t)}>0 \\
\label{def_Dbarmax}
\bar{D}_{\max}(t) &= D_{\max}-D^* + \frac{\dot{y}_\text{ref}(t)}{y_\text{ref}(t)}>0.
\end{align}
Notice that the following inequalities hold for all $\eta(t)\in\mathbb{R}$, $e_2(t)\in\mathbb{R}$, $\delta(t)\in\mathbb{R}$, $\gamma\in\mathbb{R}$ and all $t\geq0$
\begin{align}
\label{caseI_1}
\bigg|\text{sat}_{{\mathcal{D}(t)}} \Big(e_2(t)+\gamma\eta(t)+\gamma\delta(t) \Big)\bigg|&\leq \max\big(\bar{D}_\text{min}(t),\bar{D}_\text{max}(t)\big)\\
\label{caseI_2}
\bigg|\text{sat}_{{\mathcal{D}(t)}} \left(\frac{\gamma}{2}\eta(t) \right)\bigg|&\leq \max\big(\bar{D}_\text{min}(t),\bar{D}_\text{max}(t)\big).
\end{align}
With~\eqref{ValidYref} and the definitions~\eqref{def_Dbarmax},~\eqref{def_Dbarmin} we are in the position to establish a more conservative bound  of the right-hand side of \eqref{caseI_1} and \eqref{caseI_2} for all~$t\geq0$, namely 
\begin{align}
\max\big(\bar{D}_\text{min}(t),\bar{D}_\text{max}(t)\big)< D_{\max}-D_\text{min}.
\end{align}
The time derivative of the squared state~$\eta(t)$ is governed by
\begin{align}
\label{eta2_dot}
\frac{\text{d} \eta^2(t)}{\text{d}t} = -2\eta(t) \text{sat}_{\mathcal{D}(t)} \Big(e_2(t)+\gamma\eta(t)+\gamma\delta(t) \Big).
\end{align}
as a consequence of \eqref{dot_eta(t)}. We next distinguish the following cases in order to derive a bound for~\eqref{eta2_dot}:
\begin{align}
\label{def_caseI}
&\text{(I)}:&     |e_2(t)+\gamma\delta(t)|&>\frac{\gamma}{2}|\eta(t)|,&\\
\label{def_caseIIa}
&\text{(IIa)}:&  |e_2(t)+\gamma\delta(t)|&\leq\frac{\gamma}{2}|\eta(t)|&\text{and}~ \eta(t)\geq0&,\\
\label{def_caseIIb}
&\text{(IIb)}:& |e_2(t)+\gamma\delta(t)|&\leq\frac{\gamma}{2}|\eta(t)|&\text{and}~ \eta(t)<0&
\end{align}
In case (I) it follows from \eqref{caseI_2}--\eqref{def_caseI}
\begin{align}
\frac{\text{d} \eta^2(t)}{\text{d}t}+2\eta(t)\text{sat}_{{\mathcal{D}(t)}}\left(\frac{\gamma}{2}\eta(t)\right) &\leq 4|\eta(t)|  \left(D_{\max}-D_\text{min} \right)\nonumber\\
&\leq\mu_2 \big|e_2(t)+\gamma\delta(t)\big|
\end{align}
where $\mu_2 = 8\gamma^{-1} \left(D_{\max}-D_\text{min} \right)>0$. Hence, it follows from the previous inequality that
\begin{align}
\label{caseI}
\frac{\text{d} \eta^2(t)}{\text{d}t}\leq-2\eta(t)\text{sat}_{{\mathcal{D}}(t)}\left(\frac{\gamma}{2}\eta(t)\right)+\mu_2 \big|e_2(t)+\gamma\delta(t)\big|.
\end{align}
holds. Next using~\eqref{def_caseIIa} in case (IIa), the following chain of inequalities holds
\begin{align}
\label{case_IIa}
\frac{\gamma}{2}\eta(t)\leq e_2(t)+\gamma\delta(t)+\gamma\eta(t)\leq \frac{3\gamma}{2}\eta(t).
\end{align}
The fact that the saturation function in~\eqref{eta2_dot} is non-decreasing and the previous inequality imply
\begin{align}
\frac{\text{d} \eta^2(t)}{\text{d}t}\leq-2\eta(t)\text{sat}_{{\mathcal{D}(t)}}\left(\frac{\gamma}{2}\eta(t)\right),
\end{align}
such that~\eqref{caseI} also holds in case (IIa). For case~(IIb), a similar argument works when using the inequalities
\begin{align}
\label{eta2_caseIIa}
\frac{\gamma}{2}\eta(t)\geq e_2(t)+\gamma\delta(t)+\gamma\eta(t)\geq \frac{3\gamma}{2}\eta(t),
\end{align}
which imply that~\eqref{eta2_caseIIa} as well as~\eqref{caseI} also hold in the present case. Comparing all thee cases~\eqref{def_caseI}--\eqref{def_caseIIb}, we conclude that~\eqref{caseI} is a valid bound for~\eqref{eta2_dot} for all $t\geq0$. By taking~\eqref{ineq_xsatx} into account we get
\begin{align}
\label{xsatx_final}
2\eta(t)\text{sat}_{{\mathcal{D}}(t)}\left(\frac{\gamma}{2}\eta(t)\right) \geq \mu_1(y_\text{ref})\frac{\eta^2(t)}{1+\sqrt{\eta^2(t)}}
\end{align}
for all $t\geq0$ where
\begin{align}
\label{def_mu1}
\mu_1(y_\text{ref})=&\min(2,\gamma)\cdot\min\bigg\{1,D^*-D_{\min} - \sup_{\tau\geq0} \frac{\dot{y}_\text{ref}(\tau)}{y_\text{ref}(\tau)}, \nonumber\\&D_{\max}-D^* +\inf_{\tau\geq0} \frac{\dot{y}_\text{ref}(\tau)}{y_\text{ref}(\tau)}\bigg\}.
\end{align}
Eventually, \eqref{estimate_eta2_dot} straightforwardly follows from \eqref{caseI} in conjunction with \eqref{xsatx_final} and definition~\eqref{def_mu1}.$\hfill\square$
%
\bibliographystyle{plain}             
\bibliography{references}           
\parpic{\includegraphics[width=1in,clip,keepaspectratio]{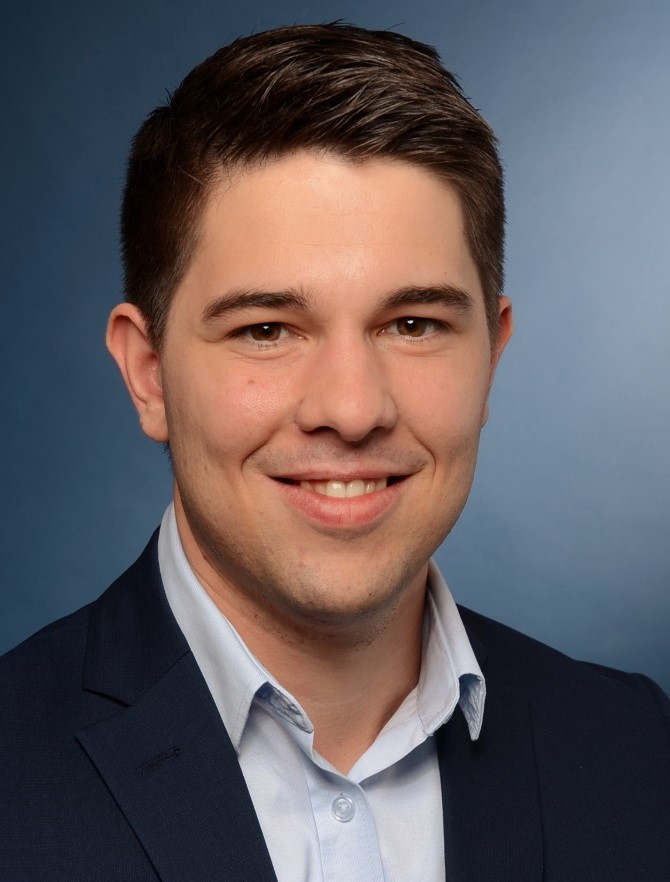}} \small\noindent {\bf Kevin Schmidt} received the B.Sc. and M.Sc. joint degrees in medical engineering from the Universities of Stuttgart and Tuebingen, Germany, in 2014 and 2016 respectively. He was with the Cymer Center for Control Systems and Dynamics, UC San Diego, USA in 2016. Since 2016, he has been a Research Assistant with the Institute for System Dynamics, University of Stuttgart, Germany. His research interests include control of distributed parameters systems and optimal control with applications to transport systems, thermomechanics and adaptive optics devices.\\[-3ex]
\parpic{\includegraphics[width=1in,clip,keepaspectratio]{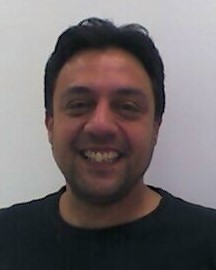}} \small \noindent {\bf Iasson Karafyllis} received a Ph.D. in Mathematics from the National Technical University of Athens (NTUA), Greece in 2003. He is Assistant Professor in the Department of Mathematics, NTUA, Greece. He is a co-author (with Zhong-Ping Jiang) of the book Stability and Stabilization of Nonlinear Systems, Springer-Verlag London (Series: Communications and Control Engineering), 2011. Since 2013 he is an Associate Editor for the International Journal of Control and for the IMA Journal of Mathematical Control and Information. His research interests include mathematical control theory and nonlinear systems theory.\\[-3ex]
\parpic{\includegraphics[width=1in,clip,keepaspectratio]{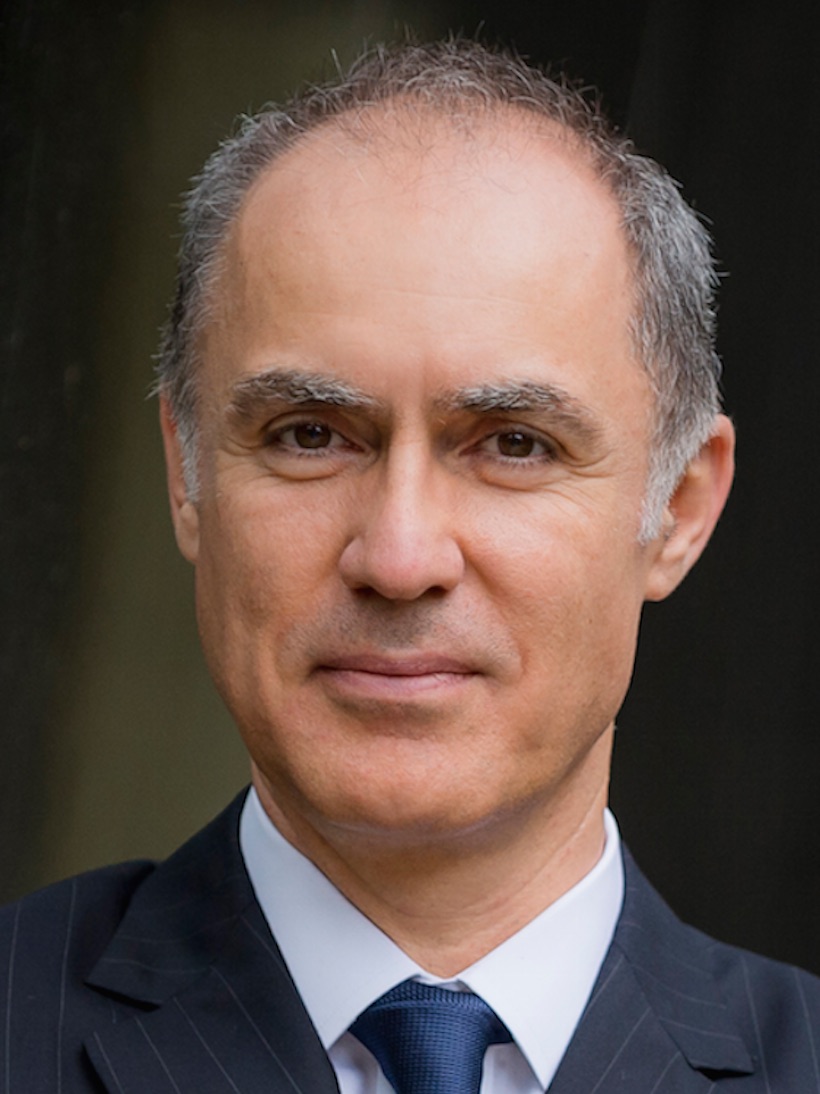}} \small \noindent {\bf Miroslav Krstic} is Fellow of IEEE, IFAC, ASME, SIAM, and IET (UK), Associate Fellow of AIAA, and foreign member of the Academy of Engineering of Serbia. He has received ASME Oldenburger Medal, ASME Nyquist Lecture Prize, ASME Paynter Outstanding Investigator Award, the PECASE, NSF Career, and ONR Young Investigator awards, the Axelby and Schuck paper prizes, the Chestnut textbook prize, and the first UCSD Research Award given to an engineer. Krstic has coauthored twelve books on adaptive, nonlinear, and stochastic control, extremum seeking, control of PDE systems including turbulent flows, and control of delay systems.
%

\end{document}